\documentclass[12pt]{article}

\usepackage{amssymb}
\usepackage{amsthm}
\usepackage{amsmath}

\def \tilde{\widetilde}

\newcommand{\st}[1]{\ensuremath{^{\scriptstyle \textrm{#1}}}}

\newcommand{\CC}{{\mathbb C}}

\newcommand{\ZZ}{{\mathbb Z}}

\newcommand{\A}{\mathcal{A}}

\renewcommand{\L}{\mathcal{L}}
\renewcommand{\O}{\mathcal{O}}

\newcommand{\fa}{{\mathfrak a}}

\newcommand{\fg}{{\mathfrak g}}

\newcommand{\fp}{{\mathfrak p}}

\newcommand{\ad}{\mathop{\rm ad}}
\newcommand{\codim}{\mathop{\rm codim}}
\newcommand{\const}{\mathop{\rm const}}

\newcommand{\Der}{\mathop{\rm Der}}
\renewcommand{\div}{\mathop{\rm div}}

\newcommand{\even}{\mathop{\rm even}}
\newcommand{\Ind}{\mathop{\rm Ind}}
\newcommand{\odd}{\mathop{\rm odd}}

\newcommand{\Vect}{\mathop{\rm Vect}}
\newcommand{\Vir}{\mathop{\rm Vir}}

\newcounter{bean}
\newenvironment{deflist}[1]%
    {
      \begin{list}{\bf #1\arabic{bean}}
         {\usecounter{bean}
              \setcounter{bean}{-1}
          \labelsep=1em
          \settowidth{\labelwidth}{#1\thebean:}
          \addtolength{\labelwidth}{1.1ex}
          \leftmargin=\labelwidth
          \addtolength{\leftmargin}{\labelsep} }

    }    {\end{list}}

\makeatletter
\renewcommand\section{\@startsection {section}{1}{\z@}%
                                   {-3.5ex \@plus -1ex \@minus -.2ex}%
                                   {2.3ex \@plus.2ex}%
                                   {\normalfont\large\bfseries}}

\@addtoreset{equation}{section}
\makeatother

\newcounter{saveeq}
{ \refstepcounter{equation}
  \setcounter{saveeq}{\value{equation}}
  \setcounter{equation}{0}
  }%
{ \setcounter{equation}{\value{saveeq}}%
  }

\newtheorem{theorem}{Theorem}%%%[section]
\newtheorem*{theorem*}{Theorem}

\newtheorem{proposition}{Proposition}
\newtheorem*{proposition*}{Proposition}

\newtheorem*{corollary*}{Corollary}
\theoremstyle{definition}
\newtheorem*{definition*}{Definition}
\newtheorem*{comments*}{Comments}

\newtheorem*{example*}{Example}

\theoremstyle{remark}
\newtheorem*{remark*}{Remark}
\newtheorem*{remarks*}{Remarks}

\newcommand{\alphaparenlist}{% changes enumerate 1st level to (a)...(z)
  \renewcommand{\theenumi}{\alph{enumi}}%
  \renewcommand{\labelenumi}{(\theenumi)}%
}
\newcommand{\Alphalist}{% changes enumerate 1st level to A. ... Z.
  \renewcommand{\theenumi}{\Alph{enumi}}%
  \renewcommand{\labelenumi}{\theenumi.}%
}
\newcommand{\arabiclist}{% changes enumerate 1st level to 1. ... 9.
  \renewcommand{\theenumi}{\arabic{enumi}}%
  \renewcommand{\labelenumi}{\theenumi.}%
}
\newcommand{\arabicparenlist}{% changes enumerate 1st level to (1) ... (9)
  \renewcommand{\theenumi}{\arabic{enumi}}%
  \renewcommand{\labelenumi}{(\theenumi)}%
}
\newcommand{\romanparenlist}{% changes enumerate 1st level to (i)...(x)
  \renewcommand{\theenumi}{\roman{enumi}}%
  \renewcommand{\labelenumi}{(\theenumi)}%
}
\newcommand{\alphaparenlistii}{% changes enumerate 2nd level to (a)...(z)
  \renewcommand{\theenumii}{\alph{enumii}}%
  \renewcommand{\labelenumii}{(\theenumii)}%
}
\newcommand{\romanparenlistii}{% changes enumerate 2nd level to (i)...(x)
  \renewcommand{\theenumii}{\roman{enumii}}%
  \renewcommand{\labelenumii}{(\theenumii)}%
}
\makeatletter
\def\@maketitle{\newpage
 \null
 \vskip 2em
 \begin{center}%
  {\@date}%
  \vskip 3em
  {\Large\bf \@title \par}%
  \vskip 1.5em
  {\normalsize
   \lineskip .5em
   \begin{tabular}[t]{c}\@author
   \end{tabular}\par}%
  \vskip 2em

 \end{center}%
 \par
 \vskip 2.5em}
\makeatother
\begin{document}
\date{\hspace{2.8in} Talk given at the conference\\
\hspace{2.9in} ``Visions in Mathematics toward \\
\hspace{2in}    the year 2000'',\\
\hspace{2.5in}  August 1999, Tel-Aviv}

\title{Classification of infinite-dimensional simple groups of
  supersymmetries and quantum field theory}

\author{Victor G. Kac\thanks{%
    Department of Mathematics, M.I.T.,
    Cambridge, MA 02139,
    $<$kac@math.mit.edu$>$
 %%     Supported in part by NSF grant DMS-9622870.
}}

\maketitle

%\begin{abstract}

%\end{abstract}

\section*{Introduction}

This work was motivated by two seemingly unrelated problems:

\begin{enumerate}
\item %%1
  Lie's problem of classification of ``local continuous
  transformation groups of a finite-dimensional manifold''.\vspace{-1ex}

\item %%2
  The problem of classification of operator product expansions
  (OPE) of chiral fields in $2$-dimensional conformal field theory.
\end{enumerate}

I shall briefly explain in \S~\ref{sec:5} how these problems are related to
each other via the theory of conformal algebras \cite{DK},
\cite{K4}--\cite{K6}.  This connection led to the classification of finite
systems of chiral bosonic fields such that in their OPE only
linear combinations of these fields and their derivatives occur
\cite{DK}, which is, basically, what is called a ``finite
conformal algebra''.

It is, of course, well known that a solution to Lie's problem
requires quite different methods in the cases of finite- and
infinite-dimensional groups.  The most important advance in the
finite-dimensional case was made by W.~Killing and E.~Cartan at
the end of the 19\st{th} century who gave the celebrated
classification of simple finite-dimensional Lie algebras over
$\CC$.  The infinite-dimensional case was studied by E.~Cartan in
a series of papers written in the beginning of the 20\st{th}
century, which culminated in his classification of
infinite-dimensional ``primitive'' Lie algebras \cite{C}.

The advent of supersymmetry in theoretical physics in the 1970's
motivated the work on the ``super'' extension of Lie's problem.
In the finite-dimensional case the latter problem was settled in
\cite{K2}.  However, it took another 20~years before the problem
was solved in the infinite-dimensional case \cite{K7},
\cite{CK2}, \cite{CK3}.  An entertaining account of the historical
background of the four classifications mentioned above may be
found in the review \cite{St}.

A large part of my talk (\S\S~\ref{sec:1}--\ref{sec:4}) is devoted
to the explanation of the fourth classification, that of simple
infinite-dimensional local supergroups of transformations of a
finite-dimensional supermanifold.  The application of this result
to the second problem, that of classification of OPE when
fermionic fields are allowed as well, or, equivalently, of finite
conformal superalgebras, is explained in \S~\ref{sec:5}.

I am convinced, however, that the classification of
infinite-dimensional supergroups may have applications to
``real'' physics as well.  The main reason for this belief is the
occurrence in my classification of certain exceptional
infinite-dimensional Lie supergroups that are natural extension of
the compact Lie groups $SU_3 \times SU_2 \times U_1$ and $SU_5$.
In \S~\ref{sec:7} I formulate a system of
axioms, which, via representation theory of  the corresponding
Lie superalgebras (see~\S~\ref{sec:6}),
produce precisely all the multiplets of fundamental particles of
the Standard model.

I wish to thank K.~Gawedzki, D.~Freedman, L.~Michel, L.~Okun, R.~Penrose,
S.~Petcov, A.~Smilga, H.~Stremnitzer, I.~Todorov for discussions
and correspondence.  I am grateful to
P.~Littlemann and J.~van der Jeugt for providing tables of certain
branching rules.

\section{Lie's problem and Cartan's theorem}
\label{sec:1}

In his ``Transformation groups'' paper \cite{L} published in
1880, Lie argues as follows.  Let $G$ be a ``local continuous
group of transformations'' of a finite-dimensional manifold.  The
manifold decomposes into a union of orbits, so we should first
study how the group acts on an orbit and worry later how the
orbits are put together.  In other words, we should first study
transitive actions of $G$.  Furthermore, even for a transitive
action it may happen that $G$ leaves invariant a fibration by
permutting the fibers.  But then we should first study how $G$
acts on fibers and on the quotient manifold and worry later how
to put these actions together.  We thus arrive at the problem of
classifications of transitive primitive (i.e.,~leaving no
invariant fibrations) actions.

Next Lie establishes his famous theorems that relate the action
of $G$ on a manifold $M$ in a neighborhood of a point $p$ to the
Lie algebra of vector fields in this neighborhood generated by
this action.  Actually, he talks about formal vector fields in a
formal neighborhood of $p$, hence $G$ gives rise to a Lie algebra
$L$ of formal vector fields and its canonical filtration by
subalgebras $L_j$ of $L$ consisting of vector fields that vanish
at $p$ up to the order $j+1$.  Then transitivity
of the action of $G$ is equivalent to the
property that $\dim L/L_0=m$:= $\dim M$, and primitivity
is equivalent to the property that
$L_0$ is a maximal subalgebra of $L$.

The first basic example is the Lie algebra of all formal vector
fields in $m$ indeterminates:
\begin{displaymath}
  W_m = \left\{ \sum^m_{i=1} P_i (x)
        \frac{\partial}{\partial x_i} \right\}
\end{displaymath}
($P_i(x)$ are formal power series in $x=(x_1 , \ldots
,x_m)$), endowed with the formal topology.  A subalgebra $L$ of $W_m$
is called transitive if $\dim L/L_0=m$, where $L_0=(W_m)_0 \cap L$.
A rigorous statement
of Lie's problem is as follows:

\textbf{$1$\st{st} formulation.}~~Classify all closed transitive
subalgebras $L$ of $W_m$ such that
$L_0$ is a maximal subalgebra of $L$, up to a continuous automorphism of
$W_m$.

E.~Cartan published a solution to this problem in the
infinite-dimensional case in 1909 \cite{C}.  The result is that a
complete list over $\CC$
(conjectured by Lie) is as follows $(m \geq 1)$:
\begin{displaymath}
  \begin{array}{lllll}
     &1. & W_m \, , \\
    &2.  &S_m & = &
    \{ X \in W_m | \div X =0 \} \,   (m \geq 2) \, ,\\
     &2'. & CS_m & = &\{ X \in W_m | \div X = \const \} \,
                (m \geq 2) \, ,\\
    &3.  &H_m &= &\{ X \in W_m | X\omega_s =0 \} \,    (m=2k)\, , \\
%%  \end{array}
%%\end{displaymath}
%
\noalign{\hspace*{-1ex}
    \hbox{where $\omega_s = \sum^k_{i=1} dx_i \wedge dx_{k+i}$
      is a symplectic form,}}\\[-2ex]
%
%%\begin{eqnarray*}
%%  \begin{array}{lllll}
     &  3'. & CH_m &= &\{ X \in W_m | X\omega_s = \const \omega_s
     \} \,
        (m=2k), \\
      & 4. &  K_m &= & \{ X \in W_m | X\omega_c =
      f\omega_c \} \, (m=2k+1) \, ,
  \end{array}
\end{displaymath}
where $\omega_c = dx_m + \sum^k_{i=1} x_i dx_{k+i} $
is a contact form and $f$ is a formal power series (depending on $X$).

The work of Cartan had been virtually forgotten until the
sixties.  A resurgence of interest in this area began with the
papers \cite{SS} and \cite{GS}, which developed an adequate
language and machinery of filtered and graded Lie algebras.  The
work discussed in the present talk uses heavily also the ideas
from \cite{W}, \cite{K1} and \cite{G2}.

The transitivity of the action implies that $L_0$ contains no non-zero ideals
of $L$. The pair $(L,L_0)$ is called primitive if $L_0$ is a (proper)
maximal subalgebra of $L$ which contains no non-zero ideals of $L$.
Using the Guillemin-Sternberg realization theorem \cite{GS}, \cite{B1}, it
is easy to show \cite{G2} that the $1$\st{st} formulation of
Lie's problem is equivalent to the following, more invariant,
formulation:

\textbf{$2$\st{nd} formulation.}  Classify all primitive pairs $(L,L_0)$,
where $L$ is a linearly compact Lie
algebra and $L_0$ is
its open  subalgebra.

Recall that a topological Lie algebra is called linearly compact
if its underlying space is isomorphic to a topological product of
discretely topologized finite-dimensional vector spaces (the basic
examples are:  finite-dimensional spaces with discrete topology
and the space of formal power series in $x$ with formal topology).

Any linearly compact Lie algebra $L$ contains an open (hence of
finite codimension) subalgebra $L_0$ \cite{G1}.  Hence, if $L$ is
simple, choosing any maximal open subalgebra $L_0$, we get a
primitive pair $(L,L_0)$.  One can show that there
exists a unique such $L_0$, and this leads to the four series
$1,2,3$ and~$4$.  The remaining series $2'$ and $3'$ are not
simple, they actually are the Lie algebras of derivations of $2$
and $3$, but the choice of $L_0$ is again unique.

Using the structure results on general transitive linearly
compact Lie algebras \cite{G1}, it is not difficult to reduce, in
the infinite-dimensional case, the classification of
primitive pairs to the classification of simple linearly compact
Lie algebras (cf.~\cite{G2}).  Such a reduction is possible also
in the Lie superalgebra case, but it is much more complicated for
two reasons:
\alphaparenlist
\vspace{-1ex}
\begin{enumerate}
\item%%a
a simple linearly compact Lie superalgebra may have several
maximal open subalgebras,
\vspace{-1ex}
\item  %%b
construction of arbitrary primitive pairs in terms of
  simple primitive pairs is more complicated in the superalgebra
  case.
\end{enumerate}
\vspace{-1ex}
In the next sections I will discuss in some detail the classification
of infinite-dimensional simple linearly compact Lie
superalgebras.

\section{Statement of the main theorem}
\label{sec:2}
The ``superization'' basically amounts to adding anticommuting
indeterminates.  In other words, given an algebra (associative or
Lie) $\A$ we consider the Grassmann algebra $\A\langle n \rangle$ in $n$
anticommuting indeterminates $\xi_1 , \ldots , \xi_n$ over $\A$.
This algebra carries a canonical $\ZZ / 2 \ZZ$-gradation, called
parity, defined by letting
\begin{displaymath}
  p(\A)=\overline{0}, \quad p (\xi_i) = \overline{1},
  \quad \overline{0}, \overline{1} \in \ZZ / 2 \ZZ \, .
\end{displaymath}
For example, $\CC\langle n \rangle$ is the Grassmann algebra in $n$
indeterminates over $\CC$.  If $\O_m$ denotes the algebra of
formal power series over $\CC$ in $m$ indeterminates, then
$\O_m\langle n \rangle$ is the algebra over $\CC$ of formal power series in $m$
commuting indeterminates $x=(x_1, \ldots , x_m)$ and $n$
anticommuting indeterminates $\xi = (\xi_1 , \ldots , \xi_n)$:
\begin{displaymath}
  x_ix_j=x_jx_i, \quad x_i\xi_j = \xi_j x_i, \quad
  \xi_i \xi_j =-\xi_j \xi_i \, .
\end{displaymath}

Recall that a derivation $D$ of parity $p(D) \in \ZZ /2\ZZ$ of a
$\ZZ /2\ZZ$-graded algebra is a vector space endomorphism
satisfying condition
\begin{displaymath}
  D(ab) = (Da)b + (-1)^{p(D)p(a)}a(Db) \, .
\end{displaymath}
Furthermore the sum of the spaces of derivations of parity
$\overline{0}$ and $\overline{1}$ is closed under the ``super''
bracket:
\begin{displaymath}
  [D,D_1] = DD_1 - (-1)^{p(D)p(D_1)} D_1D \, .
\end{displaymath}
This ``super'' bracket satisfies ``super'' analogs of
anticommutativity and Jacobi identity, hence defines what is
called a Lie superalgebra.

For example, the algebra $\A\langle n \rangle$ has derivations
$\displaystyle{\frac{\partial}{\partial \xi_i}}$ of parity~$\overline{1}$
defined by
\begin{displaymath}
  \frac{\partial}{\partial \xi_i} (a) =0 \hbox{ for }
  a \in \A, \quad \frac{\partial }{\partial \xi_i}
  (\xi_j) = \delta_{ij} \, ,
\end{displaymath}
and these derivations anticommute, so that $\displaystyle{\left[
    \frac{\partial}{\partial \xi_i} , \frac{\partial}{\partial
      \xi_j}\right]=0}$.

The ``super'' analog of the Lie algebra $W_m$ is the Lie
superalgebra, denoted by $W(m|n)$, of all continuous derivations
of the $\ZZ /2\ZZ$-graded algebra $\O_m\langle n \rangle, n \in \ZZ_+$, with
the defined above ``super''bracket:
\begin{displaymath}
  W(m|n) = \left\{ \sum^m_{i=1} P_i (x,\xi)
    \frac{\partial}{\partial x_i} + \sum^n_{j=1}
    Q_j (x,\xi) \frac{\partial}{\partial \xi_j} \right\} \, ,
\end{displaymath}
where $P_i (x,\xi), Q_j (x,\xi) \in \O_m\langle n \rangle$.  In a more
geometric language, this is the Lie superalgebra of all formal
vector fields on a supermanifold of dimension $(m|n)$.

There is a unique way to extend divergence from $W_m$ to
$W(m|n)$ such that the divergenceless vector fields form a
subalgebra:
\begin{displaymath}
  \div \left(\sum_i P_i \frac{\partial}{\partial x_i} + \sum_j
    Q_j \frac{\partial}{\partial \xi_j} \right)
  = \sum_j \frac{\partial P_i}{\partial x_i} + \sum_j (-1)^{p(Q_j)}
  \frac{\partial Q_j}{\partial \xi_j} \, ,
\end{displaymath}
and the ``super'' analog of $S_m$ is
\begin{displaymath}
  S(m|n) = \{ X \in W(m|n) |  \div X =0 \} \, .
\end{displaymath}

In order to define ``super'' analogs of the Hamiltonian and
contact Lie algebras $H_m$ and $K_m$, introduce a ``super'' analog
of the algebra of differential forms \cite{K2}.  This is an
associative algebra over $\O_m\langle n \rangle$, denoted by $\Omega(m|n)$,
on generators $dx_1 , \ldots , dx_m$, $d\xi_1,\ldots,d\xi_n$
and defining relations:
\begin{displaymath}
  dx_idx_j =-dx_jdx_i, \quad dx_i d\xi_j=d\xi_j dx_i, \quad
  d\xi_i d\xi_j = d\xi_j d\xi_i \, ,
\end{displaymath}
and the $\ZZ /2\ZZ$ gradation defined by:
\begin{displaymath}
  p(x_i)=p(d\xi_j)=\overline{0}, \quad
  p(\xi_j)=p(dx_i)=\overline{1} \, .
\end{displaymath}
The algebra $\Omega(m|n)$ carries a unique continuous derivation
$d$ of parity $\overline{1}$ such that
\begin{displaymath}
  d(x_i)=dx_i, \quad d(\xi_j)=d\xi_j, \quad d(dx_i)=0,\quad
  d(d\xi_j)=0 \, .
\end{displaymath}
The operator $d$ has all the usual properties, e.g.:
\begin{displaymath}
  df= (-1)^{p(f)}\sum_i \frac{\partial f}{\partial x_i} \, dx_i +
  \sum_j \frac{\partial f}{\partial \xi_j} \, d\xi_j
  \hbox{ for }f \in \O_m\langle n \rangle, \hbox{ and } d^2=0\, .
\end{displaymath}
As usual, for any $X \in W(m|n)$ one defines a derivation
$\iota_X$ (contraction along $X$) of the algebra $\Omega(m|n)$
by the properties (here $x$ stands for $x$ and $\xi$):
\begin{displaymath}
  p(\iota_X) = p(X) + \overline{1}, \quad
  \iota_X (x_j)=0, \quad
  \iota_X (dx_j) =(-1)^{p(X)} X(x_j) \, .
\end{displaymath}
The action of any $X \in W(m|n)$ on $\O_m\langle n \rangle$ extends in a
unique way to the action by a derivation of $\Omega(m|n)$ such
that $[X,d]=0$.  This is called Lie's derivative and is usually denoted
by $L_X$, but we shall write $X$ in place of $L_X$ unless
confusion may arise.  One has the usual Cartan's formula for this
action:  $L_X = [d,\iota_X]$.

Using this action, one can define super-analogs of the
Hamiltonian and contact Lie algebras for any $n \in \ZZ_+$:
\begin{eqnarray*}
 H(m|n) &=& \{ X \in W(m|n) |  X\omega_s =0 \} \, , \\
\noalign{\hbox{\hbox{where } $\omega_s = \sum^k_{i=1} dx_i \wedge dx_{k+i}
      + \sum^n_{j=1} (d \xi_j)^2$,}}\\[-2ex]
  K(m|n) &=& \{ X \in W(m|n)|  X \omega_c = f\omega_c \} \, , \\
 \noalign{\hbox{\hbox{where }$\omega_c = dx_m + \sum^k_{i=1} x_i dx_{k+i}
      + \sum^n_{j=1} \xi_j \, d\xi_j ,  \hbox{ and } f\in  \O_m\langle
n \rangle$.}}
\end{eqnarray*}

Note that $W(0|n)$, $S(0|n)$ and $H(0|n)$ are
finite-dimensional Lie superalgebras.  The Lie superalgebras
$W(0|n)$ and $S(0|n)$ are simple iff $n \geq 2$ and $n \geq 3$,
respectively.  However, $H(0|n)$ is not simple as its derived
algebra $H'(0|n)$ has codimension~$1$ in $H(0|n)$, but
$H'(0|n)$ is simple iff $n\geq 4$.  Thus, in the Lie
superalgebra case the lists of simple finite- and
infinite-dimensional algebras are much closer related than in the
Lie algebra case.

The four series of Lie superalgebras are infinite-dimensional if
$m \geq 1$, in which case they are simple except for $S(1|n)$.
The derived algebra $S'(1|n)$ has codimension~$1$ in $S(1|n)$,
and $S'(1|n)$ is simple iff $n \geq 2$.

In my paper \cite{K2} I conjectured that the constructed above
four series exhaust all infinite-dimensional simple linearly
compact Lie superalgebras.  Remarkably, the situation turned out
to be much more exciting.

As was pointed out by several mathematicians, the Schouten
bracket \cite{Sc} makes the space of polyvector fields on a
$m$-dimensional manifold into a Lie superalgebra.  The formal
analog of this is the following fifth series of superalgebras,
called by physicists the Batalin-Vilkoviski algebra:
\begin{displaymath}
  HO(m|m) =\{ X \in W(m|m) |   X\omega_{os} =0 \} \, ,
\end{displaymath}
where $\omega_{os} = \sum^m_{i=1} \, dx_i d\xi_i$ is an odd
symplectic form.  Furthermore, unlike in the $H(m|n)$ case, not
all vector fields of $HO(m|n)$ have zero divergence, which gives
rise to the sixth series:
\begin{displaymath}
  SHO(m|m) =\{ X \in HO(m|m)|  \div X=0 \} \, .
\end{displaymath}
The seventh series is the odd analog of $K(m|n)$ \cite{ALS}:
\begin{displaymath}
KO(m|m+1) =\{ X \in W(m|m+1) |  X\omega_{oc} =f\omega_{oc} \} \, ,
\end{displaymath}
where $\omega_{oc}=d\xi_{m+1} + \sum^m_{i=1} (\xi_i \, dx_i + x_i \, d
\xi_i)$ is an odd contact form.  One can take again the divergence
$0$ vector fields in $KO(m|m+1)$ in order to construct the
eighth series, but the situation is more interesting.
It turns out that for each $\beta \in \CC$ one can define the
deformed divergence $\div_{\beta}X$ \cite{Ko}, \cite{K7}, so that $\div = \div_0$ and
\begin{displaymath}
  SKO(m|m+1;\beta) =\{ X \in KO(m|m+1) |   \hbox{div}_{\beta} X=0 \}
\end{displaymath}
is a subalgebra.  The superalgebras $HO(m|m)$ and $KO(m|m+1)$
are simple iff $m \geq 2$ and $m \geq 1$, respectively.  The
derived algebra $SHO'(m|m)$ has codimension~$1$ in $SHO(m|m)$,
and it is simple iff $m \geq 3$.  The derived algebra
$SKO'(m|m+1;\beta)$ is simple iff $m \geq 2$, and it coincides
with $SKO(m|m+1;\beta)$ unless $\beta =1$ or $\frac{m-2}{m}$
when it has codimension~$1$.

Some of the examples described above have simple ``filtered
deformations'', all of which can be obtained by the following
simple construction.  Let $L$ be a subalgebra of $W(m|n)$, where
$n$ is even.  Then it happens in three cases that
\begin{displaymath}
  L^{\sim}:=(1+\Pi^n_{j=1} \xi_j)L
\end{displaymath}
is different from $L$, but is closed under bracket.  As a result
we get the following three series of superalgebras:
$S^{\sim}(0|n)$ \cite{K2},
$SHO^{\sim}(m|m)$ \cite{CK2} and
$SKO^{\sim}(m|m+1;\frac{m+2}{m})$ \cite{Ko} (the constructions in
\cite{Ko} and \cite{CK2} were more complicated).  We thus get the
ninth and the tenth series of simple infinite-dimensional Lie
superalgebras:
\begin{eqnarray*}
  SHO^{\sim}(m|m), \quad m \geq 2, \,\, m \even \, ,\\
SKO^{\sim}(m|m+1;\frac{m+2}{m}),\,\, m \geq 3, m \odd \, .
\end{eqnarray*}

It is appropriate to mention here that the four series $W(0|n)$,
$S(0|n)$, $S^{\sim}(0|n)$ and $H'(0|n)$ along with the
classical series $s\ell(m|n)$ and $osp(m|n)$, strange series
$p(n)$ and $q(n)$, two exceptional superalgebras of dimension $40$
and $31$ and a family of $17$-dimensional exceptional
superalgebras along with the marvelous five exceptional Lie
algebras, comprise a complete list of simple finite-dimensional
Lie superalgebras \cite{K2}.

A surprising discovery was made in \cite{Sh1}
where the existence of three exceptional simple
infinite-dimensional Lie superalgebras was announced.  The proof
of the existence along with one more exceptional example was
given in \cite{Sh2}.  An explicit construction of these four
examples was given later in \cite{CK3}.  The fifth exceptional
example was found in the work on conformal algebras \cite{CK1}
and independently in \cite{Sh2}.  (The alleged sixth exceptional
example $E(2|2)$
of \cite{K7} turned out to be isomorphic to $SK0(2|3;1)$ \cite{CK3}.)

Now I can state the main theorem.

\begin{theorem} \cite{K7}
  \label{th:1}
  The complete list of simple infinite-dimensional linearly
  compact Lie superalgebras consists of ten series of examples
  described above and five exceptional examples:  $E(1|6)$,
  $E(3|6)$, $E(3|8)$, $E(4|4)$, and $E(5|10)$.
\end{theorem}

It happens that all infinite-dimensional simple linearly compact
Lie algebras $L$ have a unique transitive primitive action \cite{G2}.  This
is certainly false in the Lie superalgebra case.  However, if $L$
is a simple linearly compact Lie superalgebra of type $X(m|n)$,
it happens that $m$ is minimal such that $L$ acts on a
super-manifold of dimension $(m|n)$ (i.e.,~$L \subset W(m|n)$),
$n$ is minimal for this $m$, and $L$ has a unique action with
such minimal $(m|n)$.  Incidentally, in all cases the growth of
$L$ equals $m$. (Recall that growth is the minimal $m$ for which
$\dim L/L_j$ is bounded by $P(j)$, where $P$ is a polynomial of
degree~$m$.)

Let me now describe  those linearly compact infinite-dimensional
Lie superalgebras $L$ that allow a transitive primitive action.
Let $S$ be a simple linearly compact infinite-dimensional Lie
superalgebra and let $S\langle n \rangle$ denote, as before, the Grassmann
algebra over $S$ with $n$ indeterminates.  The Lie superalgebra $
\Der (S\langle n \rangle)$ of all derivations of the Lie superalgebra $S\langle n \rangle$
is the following semi-direct sum:
\begin{displaymath}
  \Der (S\langle n \rangle) = (\Der S)\langle n \rangle +  W(0|n) \, .
\end{displaymath}
(For a description of $\Der S$ see \cite{K7}, Proposition~6.1.)  Denote by $\L (S,n)$ the set of all open subalgebras $L$ of $\Der
(S\langle n \rangle)$ that contain $S\langle n \rangle$ and have the property that the
canonical image of $L$ in $W(0|n)$ is a
transitive subalgebra.

Using  a description of semi-simple linearly compact Lie
superalgebras similar to the one given by Theorem~6 from
\cite{K2} (cf.~\cite{Ch}, \cite{G1} and \cite{B2}) and Proposition~4.1 from
\cite{G1}, it is easy to derive the following result.

\begin{proposition}
  \label{prop:1}
If a linearly compact infinite-dimensional Lie superalgebra $L$
allows a transitive primitive action, then $L$ is one of the
algebras of the sets $\L (S,n)$.
\end{proposition}

\begin{example*}
%  \label{ex:1}
Consider the semidirect sum $L=S\langle n \rangle + R$, where $R$ is a transitive
subalgebra of $W(0|n)$.  Then ($L$,$L_0
=S_0\langle n \rangle+ R$) is a primitive pair
if $S_0$ is a maximal open subalgebra of $S$, and these are
all primitive pairs in the case when $S=\Der S$.  One can also
replace in this construction $S$ by $\Der S$ and $S_0$ by a maximal
open subalgebra of $\Der S$ having no non-zero ideals of $\Der S$.
\end{example*}

\section{Explanation of the proof of Theorem~\ref{th:1}}
\label{sec:3}

Step 1.~~Introduce Weisfeiler's filtration \cite{W} of $L$.  For
  that choose a maximal open subalgebra $L_0$ of $L$ and a
  minimal subspace $L_{-1}$ satisfying the properties:
  \begin{displaymath}
    L_{-1} \supsetneqq L_0 , \quad [L_0,L_{-1}] \subset L_{-1} \, .
  \end{displaymath}
Geometrically this corresponds to a choice of a primitive action
of $L$ and an invariant irreducible differential system.  The
pair $L_{-1},L_0$ can be included in a unique filtration:
\begin{displaymath}
  L=L_{-d} \supset L_{-d+1} \supset \cdots \supset
  L_{-1} \supset L_0 \supset L_1 \supset \cdots \, .
\end{displaymath}
Of course, if $L$ leaves invariant no non-trivial differential
system, then the ``depth'' $d=1$ and the Weisfeiler filtration
coincide with the canonical filtration.  Incidentally, in the Lie
algebra case, $d>1$ only for $K_n$ (when $d=2$), but in the Lie
superalgebra case, $d>1$ in the majority of cases.

The associated to Weisfeiler's filtration $\ZZ$-graded Lie superalgebra is of the form $GrL
= \Pi_{j \geq -d} \fg_j$, and has the following properties:

\begin{deflist}{G}
  \item %%G0
    $\dim \fg_j < \infty$ (since $\codim L_0 < \infty$),

  \item %%G1
    $\fg_{-j} = \fg^j_{-1}$ for $ j \geq 1$ (by maximality of $L_0$),

  \item %%G2
    $[x, \fg_{-1}]=0$ for $x \in \fg_j$, $j \geq 0 \Rightarrow
    x=0$ (by simplicity of $L$),

    \item %%G3
       $\fg_0$-module $\fg_{-1}$ is irreducible (by choice of $L_{-1}$).
\end{deflist}

Weisfeiler's idea was that property (G3) is so restrictive, that
it should lead to a complete classification of $\ZZ$-graded Lie
algebras satisfying (G0)--(G3).  (Incidentally, the
infinite-dimensionality of $L$ and hence of $Gr L$, since $L$ is
simple, is needed only in order to conclude that $\fg_1 \neq 0$.)
 This indeed turned out to be the case \cite{K1}.  In fact, my
idea was to replace the condition of finiteness of the depth by
finiteness of the growth, which allowed one to add to the
Lie-Cartan list some new Lie algebras, called nowadays affine
Kac-Moody algebras.

However, unlike in the Lie algebra case, it is impossible to
classify all finite-dimensional irreducible faithful
representations of Lie superalgebras.  One needed a new idea to
make this approach work.

Step 2.~~The main new idea is to choose $L_0$ to
  be invariant with respect to all inner automorphisms of $L$
  (meaning to contain all even $\ad$-exponentiable elements of
  $L$).  A non-trivial point is the existence of such $L_0$.
  This is proved by making use of the characteristic
  supervariety, which involves rather difficult arguments of
  Guillemin \cite{G2}, that, unfortunately, I was unable to
  simplify.

Next, using a normalizer trick of Guillemin \cite{G2}, I prove,
for the above choice of $L_0$, the following very powerful
restriction on the $\fg_0$-module $\fg_{-1}$:

\begin{deflist}{G}
\setcounter{bean}{3}
   \item %%G4
       $[\fg_0, x] =\fg_{-1}$ for any non-zero even element $x$ of
       $\fg_{-1}$.

\end{deflist}

Step 3.~~Consider a faithful irreducible
  representation of a Lie superalgebra $\fp$ in a
  finite-dimensional vector space $V$.  This representation is
  called strongly transitive if
  \begin{displaymath}
    \fp \cdot x = V
 \hbox{ for any non-zero even element } x\in V \, .
\end{displaymath}
Note that property (G4) along with (G0), (G2) for $j=0$ and (G3),
shows that the $\fg_0$-module $\fg_{-1}$ is strongly transitive.

In order to demonstrate the power of this restriction, consider
first the case when $\fp$ is a Lie algebra and $V$ is purely even.
Then the strong transitivity simply means that $V \backslash
\{0\}$ is a single orbit of the Lie group $P$ corresponding to $\fp$.
It is rather easy to see that the only strongly transitive
subalgebras $\fp$ of $g \ell_{V}$ are $g\ell_{V}$, $s \ell_{V}$,
$sp_{V}$ and $csp_{V}$.  These four cases lead to $Gr L$, where
$L= W_n$, $S_n$, $H_n$ and $K_n$, respectively.

In the super case the situation is much more complicated.  First
we consider the case of ``inconsistent gradation'', meaning that
$\fg_{-1}$ contains a non-zero even element.  The classification
of such strongly transitive modules is rather long and the answer
consists of a dozen series and a half dozen exceptions (see
\cite{K7}, Theorem~3.1).  Using similar restrictions on
$\fg_{-2}, \fg_{-3}, \ldots $, we obtain a complete list of
possibilities for
\begin{displaymath}
  GrL_{\leq}:= \oplus_{j \leq 0} \fg_j
\end{displaymath}
in the case when $\fg_{-1}$ contains non-zero even elements.  It
turns out that all but one exception are not exceptions at
all, but correspond to the beginning members of some series.  As
a result, only $E(4|4)$ ``survives'' (but the infamous $E(2|2)$
doesn't).

Step 4.~~Next, we turn to the case of a consistent
  gradation, i.e.,~when $\fg_{-1}$ is purely odd.  But then
  $\fg_0$ is an ``honest'' Lie algebra, having a faithful
  irreducible representation in $\fg_{-1}$ (condition (G4)
  becomes vacuous).  An explicit description of such
  representations is given by the classical Cartan-Jacobson
  theorem.  In this case I use the ``growth'' method developed in
  \cite{K1}  and \cite{K2} to determine a complete list of
  possibilities for $Gr L_{\leq}$.  This case produces mainly the
  (remaining four) exceptions.

Step 5~is rather long and tedious \cite{CK3}  .
  For each $GrL_{\leq}$ obtained in Steps~3 and~4 we determine
  all possible ``prolongations'', i.e.,~infinite-dimensional
  $\ZZ$-graded Lie superalgebras satisfying (G2), whose negative part
  is the given $Gr L_{\leq}$.

Step 6.~~It remains to reconstruct $L$ from $GrL$,
  i.e.,~to find all possible filtered simple linearly compact Lie
  superalgebras $L$ with given $Gr L$ (such an $L$ is called a
  simple filtered deformation of $GrL$).  Of course,
  there is a trivial filtered deformation:  $GrL :=
  \Pi_{j \geq -d} \fg_j$, which is simple iff $GrL$ is.

It is proved in \cite{CK2} by a long and tedious calculation that
only $SHO(m|m)$ for $m \even \geq 2$ and $SKO(m|m+1;
  \frac{m+2}{m})$ for $m \odd \geq 3$ have a non-trivial simple
filtered deformation, which are the ninth and tenth series.  It
would be nice to have a more conceptual proof. Recall that
$SH0(m|m)$ is not simple, though it does have a simple filtered
deformation.  Note also that in the Lie algebra case all filtered
deformations are trivial.

\section{Construction of exceptional linearly compact Lie superalgebras}
\label{sec:4}

In order to describe the construction of the exceptional
infinite-dimensional Lie superalgebras (given in \cite{CK3}), I
need to make some remarks.  Let $\Omega_m = \Omega(m|0)$ be the
algebra of differential forms over $\O_m$, let
$\Omega^k_{m}$ denote the space of forms of degree $k$, and
$\Omega^k_{m,c\ell}$ the subspace of closed forms.  For any
$\lambda \in \CC$ the representation of $W_m$ on $\Omega^k_m$ can
be ``twisted'' by letting
\begin{displaymath}
  X \mapsto L_X + \lambda \div X , \quad X \in W_m \, ,
\end{displaymath}
to get a new $W_m$-module, denoted by $\Omega^k_m(\lambda)$ (the
same can be done for $W_{m,n}$).  Obviously,
$\Omega^k_m(\lambda)=\Omega^k_m$ when restricted to $S_m$.  Then
we have the following obvious $W_m$-module isomorphisms:
$\Omega^0_m \simeq \Omega^m_m (-1)$ and $\Omega^0_m(1) \simeq
\Omega^m_m$.  Furthermore, the map $X \mapsto \iota_X (dx_1
\wedge \ldots \wedge dx_m)$ gives the following $W_m$-module and
$S_m$-module isomorphisms:
\begin{displaymath}
  W_m \simeq \Omega^{m-1}_m(-1), \quad S_m \simeq \Omega^{m-1}_{m,c\ell} \, .
\end{displaymath}
We shall identify the representation spaces via these
isomorphisms.

The simplest is the construction of the largest
exceptional Lie superalgebra $E(5|10)$.  Its even part is the
Lie algebra $S_5$, its odd part is the space of closed $2$-forms
$\Omega^2_{5,c\ell}$.  The remaining commutators are defined as
follows for $X \in S_5, \quad \omega,\omega' \in \Omega^2_{5,c\ell}$:
\begin{displaymath}
  [X,\omega]=L_X\omega, \quad [\omega,\omega']=\omega \wedge \omega' \in \Omega^4_{5,c\ell}=S_5
  \, .
\end{displaymath}

Each quintuple of integers $(a_1,a_2,\ldots ,a_5)$ such that
$a=\sum_i a_i$ is even, defines a $\ZZ$-gradation of $E(5|10)$
by letting:
\begin{displaymath}
  \deg x_i=- \frac{\partial}{\partial x_i} =a_i , \quad
  \deg dx_i =a_i -\tfrac{1}{4} a \, .
\end{displaymath}
The quintuple $(2,2, \ldots ,2)$ defines the (only) consistent
$\ZZ$-gradation, which has depth $2$:    $E(5|10) =\Pi_{j \geq
  -2} \fg_j$, and one has:
\begin{displaymath}
\fg_0 \simeq s \ell_5 \hbox{ and } \fg_{-1} \simeq \Lambda^2 \CC^5,
\quad \fg_{-2} \simeq \CC^{5*} \hbox{ as $\fg_0$-modules.}
\end{displaymath}
Furthermore, $\Pi_{j \geq 0} \fg_j$ is a maximal open
subalgebra of $E(5|10)$ (the only one which is invariant with
respect to all automorphisms).  There are three other maximal open subalgebras in $E(5|10)$,
associated to $\ZZ$-gradations corresponding to quintuples
$(1,1,1,1,2)$, $(2,2,2,1,1)$ and
$(3,3,2,2,2)$, and one can show
that these four are all, up to conjugacy, maximal open
subalgebras (cf.~\cite{CK3}).

Another important $\ZZ$-gradation of $E(5|10)$, which is, unlike
the previous four, by infinite-dimensional subspaces,
corresponds to the quintuple $(0,0,0,1,1)$ and has depth~$1$:
$E(5|10) =\Pi_{\lambda \geq -1} \fg^{\lambda}$.  One has:
$\fg^0 \simeq E(3|6)$ and the $\fg^{\lambda}$ form an important
family of irreducible $E(3|6)$-modules \cite{KR}.  The
consistent $\ZZ$-gradation of $E(5|10)$ induces that of
$\fg^0:E(3|6)=\Pi_{j \geq -2} \fa_j$, where
\begin{displaymath}
  \fa_0 \simeq s\ell_3 \oplus s\ell_2 \oplus g\ell_1 , \quad
  \fa_{-1} \simeq\CC^3 \boxtimes \CC^2 \boxtimes \CC, \quad
  \fa_{-2}  \simeq \CC^3 \boxtimes \CC \boxtimes \CC \, .
\end{displaymath}
A more explicit construction of $E(3|6)$ is as follows
\cite{CK3}:  the even part is $W_3 + \Omega^0_3 \otimes
s\ell_2$, the odd part is $\Omega^1_3 (-\tfrac{1}{2}) \otimes
\CC^2$ with the obvious action of the even part, and the bracket
of two odd elements is defined as follows:
\begin{displaymath}
  [\omega \otimes u,\omega' \otimes v] =
  (\omega \wedge \omega') \otimes (u \wedge v) +
  (d\omega \wedge \omega' + \omega \wedge d\omega') \otimes (u \cdot v) \, .
\end{displaymath}
Here the identifications $\Omega^2_3 (-1) =W_3$ and
$\Omega^0_3=\Omega^3_3 (-1)$ are used.

The gradation of $E(5|10)$
corresponding to the quintuple $(0,1,1,1,1)$ has depth~$1$ and
its $0$\st{th} component is isomorphic to $E(1|6)$
(cf.~\cite{CK3}).

The construction of $E(4|4)$ is also very simple \cite{CK3}:
The even part is $W_4$, the odd part is $\Omega^1_4
(-\tfrac{1}{2})$ and the bracket of two odd elements is:
\begin{displaymath}
  [\omega,\omega'] = d\omega \wedge \omega' + \omega \wedge d\omega' \in \Omega^3_4(-1)=W_4 \, .
\end{displaymath}
The construction of $E(3|8)$ is slightly more complicated, and
we refer to \cite{CK3} for details.

\section{Classification of superconformal algebras}
\label{sec:5}

Superconformal algebras have been playing an important role in
superstring theory and in conformal field theory.  Here I will
explain how to apply Theorem~\ref{th:1} to the classification of
``linear'' superconformal algebras.  By a (``linear'')
superconformal algebra I mean a Lie superalgebra $\fg$ spanned by
coefficients of a finite family $F$ of pairwise local fields such
that the following two properties hold:

\arabicparenlist

\begin{enumerate}
\item %%1
  for $a,b \in F$ the singular part of OPE is finite, i.e.,
  \begin{displaymath}
    [a(z),b(w)] = \sum_j c_j (w) \partial^j_w \delta (z-w) \quad
    \hbox{(a finite sum),}
  \end{displaymath}
where all $c_j (w) \in \CC [\partial_w]F$ \, ,

\item %%2
  $\fg$ contains no non-trivial ideals spanned by coefficients of
  fields from a $\CC [\partial_w]$-submodule of
  $\CC [\partial_w]F$.
\end{enumerate}

This problem goes back to the physics paper \cite{RS}, some
progress in its solution was made in \cite{K6} and a complete
solution was stated in \cite{K4}, \cite{K5}.  (A complete
classification even in the ``quadratic'' case seems to be a much
harder problem, see~\cite{FL} for some very interesting
examples.)  The simplest example is the loop algebra
$\tilde{\fg}= \CC [x,x^{-1}] \otimes \fg$ (= centerless affine
Kac-Moody (super)algebra), where $\fg$ is a simple
finite-dimensional Lie (super)algebra.  Then $\displaystyle{F=\{
  a(z) = \sum_{n \in \ZZ} (x^n \otimes a)z^{-n-1}\}_{a \in
    \fg}}$, and $[a(z),b(w)]=[a,b](w) \delta (z-w)$.  The next
example is the Lie algebra $\Vect \CC^{\times}$ of regular vector
fields on $\CC^{\times}$ (= centerless Virasoro algebra);
$F$ consists of one field, the Virasoro field $  L(z)=-\sum_{n \in \ZZ} (x^n \frac{d}{dx})z^{-n-2}$, and $
         [L(z),L(w)]=\partial_w L(w) \delta (z-w)+2L(w) \delta'_w(z-w)$.

One of the main theorems of \cite{DK} states that these are all
examples in the Lie algebra case.  The strategy of the proof is
the following.  Let $\partial =\partial_z$ and consider the
(finitely generated) $\CC [\partial]$-module $R=\CC
[\partial]F$.  Define the ``$\lambda$-bracket'' $R \otimes R \to
\CC [{\lambda}] \otimes R$ by the formula:
\begin{displaymath}
  [a_{\lambda}b]=\sum_j \lambda^j c_j \, .
\end{displaymath}
This satisfies the axioms of a conformal (super)algebra (see
\cite{DK}, \cite{K4}), similar to the Lie (super)algebra axioms:

\romanparenlist
\begin{enumerate}
\item %%i
  $[\partial a_{\lambda}b] =-\lambda [a_{\lambda}b]$, $[a_{\lambda}
  \partial b] = (\partial + \lambda) [a_{\lambda}b]$,

\item %%ii
  $[a_{\lambda}b]=-[b_{-\lambda-\partial}a]$,

\item %%iii
  $[a_{\lambda}[b_{\mu}c]] = [[a_{\lambda}b]_{\lambda + \mu }c] +
  (-1)^{p(a)p(b)} [b_{\mu}[a_{\lambda}c]]$.
\end{enumerate}

The main observation of \cite{DK} is that a conformal
(super)algebra is completely determined by the Lie (super)algebra
spanned by all coefficients of negative powers of $z$ of the
fields $a(z)$, called the annihilation algebra, along with an
even surjective derivation of the annihilation algebra.
Furthermore, apart from the case of current algebras, the
completed annihilation algebra turns out to be an
infinite-dimensional simple linearly compact Lie (super)algebra of growth~$1$.
Since in the Lie algebra case the only such example is $W_1$, the
proof is finished.

In the superalgebra case the situation is much more interesting
since there are many infinite-dimensional simple linearly compact Lie
superalgebras of growth $1$.  By Theorem~\ref{th:1}, the complete
list is as follows:
\begin{displaymath}
  W(1|n), \quad S'(1|n), \quad K(1|n) \,\, \hbox{ and } E(1|6)
  \, .
\end{displaymath}
The corresponding superconformal algebras in the first three cases
are defined in the same way, except that we replace $\O_1\langle n
\rangle$ by
$\CC ((x))\langle n\rangle$; denote them by $W_{(n)}$, $S_{(n)}$ and
$K_{(n)}$, respectively.  The superconformal algebras $W_{(n)}$
and $K_{(n)}$ are simple for $n \geq 0$, except for $K_{(4)}$
which should be replaced by its derived algebra $K'_{(4)}$, and $
S'_{(n)}$ is simple for $n \geq 2$.  The unique superconformal
algebra corresponding to $E(1|6)$ is denoted by $CK_{(6)}$.  Its
construction is more difficult and may be found in \cite{CK3} or
\cite{K6}.

However, the superconformal algebra with the annihilation algebra
$S'(1|n)$ is not unique since there  are two up to conjugacy
even surjective derivations $\partial$ of $S'(1|n)$.  In order
to show this, we may assume that both $\partial$ and $S'(1|n)$
are in $W(1|n)$.  Using a change of indeterminates, we may
assume that $\displaystyle{\partial = \frac{\partial}{\partial
    x}}$, but then the standard volume form $v$ that defines
$\div$, may change to $P(\xi)v$ (since it must be annihilated by
$\frac{\partial}{\partial x}$).  Further change of indeterminates
brings this form to $(1+ \epsilon\xi_1 \ldots \xi_n)v$, where
$\epsilon =0$ or $1$.  This
gives us deformations of $S'_{(n)}$ with the annihilation
algebra $S'(1|n)$, which are derived algebras of
\begin{displaymath}
  S_{(n), \epsilon,a}=\{X\in W_{(n)} |
  \div (e^{ax}(1+\epsilon \xi_1 \ldots \xi_n)X)=0 \}\, ,
  \,\, a \in \CC \, .
\end{displaymath}
(The situation is more interesting in the case
$n=2$, since the algebra of outer derivations of $S'(1|2)$ is
$3$-dimensional \cite{P}, but this gives no new
superconformal algebras.)
One argues similarly in the case $K(1|n)$. The case $E(1|6)$ is checked
directly.  We thus have arrived at the following theorem.

\begin{theorem}
  \label{th:2}
  A complete list of superconformal algebras consists of loop
  algebras $\tilde{\fg}$, where $\fg$ is a simple
  finite-dimensional Lie superalgebra and of Lie superalgebras
  $(n \in \ZZ_+)$:  $W_{(n)}$, $S'_{(n+2),\epsilon,a}$ ($n$ even and $a=0$ if
  $\epsilon =1$),
   $K_{(n)} (n \neq 4)$, $K'_{(4)}$, and $CK_{(6)}$.

\end{theorem}

Note that the first members of the above series are well-known
superalgebras:  $W_{(0)}\simeq K_{(0)}$ is the Virasoro algebra,
$K_{(1)}$ is the Neveu-Schwarz algebra,
$K_{(2)}\simeq W_{(1)}$ is the $N=2$ algebra, $K_{(3)}$ is the
$N=3$ algebra, $S'_{(2)}$=$S'_{(2),0,0}$ is the $N=4$ algebra, $K'_{(4)}$ is
the big $N=4$ algebra.  These algebras,
along with $W_{(2)}$ and $CK_{(6)}$ are the only superconformal
algebras for which all fields are primary with positive conformal weights
\cite{K6}.  It is interesting to note that all of them are
contained in $CK_{(6)}$, which consists of $32$~fields, the even
ones are the Virasoro fields and $15$~currents that form
$\tilde{so}_6$, and the odd ones are $6$ and $10$ fields of
conformal weight $3/2$ and $1/2$, respectively.  Here is the
table of inclusions, where in square brackets the number of
fields is indicated:
\begin{eqnarray*}
  \begin{array}{cccc}
CK_{(6)} [32] & \supset & W_{(2)} [12] &
\supset W_{(1)} =K_{(2)}[4] \supset K_{(1)} [2] \supset \Vir\\
\cup && \cup\\
K_{(3)} [8] \subset K'_{(4)}[16]\qquad && S'_{(2),\epsilon,a}[8]
  \end{array}\, .
\end{eqnarray*}
All of these Lie superalgebras have a unique non-trivial central
extension, except for $K'_4$ that has three \cite{KL} and  $CK_{(6)}$
that has none.  All other
superalgebras listed by Theorem~\ref{th:2} have no non-trivial
central extensions.
(The presence of a central term is necessary for
the construction of an interesting conformal field theory.)

\section{Representations of linearly compact Lie superalgebras}
\label{sec:6}

By a representation of a linearly compact Lie superalgebra $L$ we
shall mean a continuous representation in a vector space $V$ with
discrete topology (then the contragredient representation is a
continuous representation in a linearly compact space $V^*$).
Fix an open subalgebra $L_0$ of $L$.  We shall assume that $V$ is
locally $L_0$-finite, meaning that any vector of $V$ is contained
in a finite-dimensional $L_0$-invariant subspace (this property
actually often implies that $V$ is continuous).  These kinds of
representations were studied in the Lie algebra case by Rudakov
\cite{R}.

It is easy to show that such an irreducible
$L$-module $V$ is a quotient of an induced module
$\Ind^L_{L_0}U=U(L) \otimes_{U(L_0)}U$, where $U$ is a
finite-dimensional irreducible $L_0$-module, by a (unique in good cases)
maximal submodule.  The induced module $\Ind^L_{L_0}U$ is called
degenerate if it is not irreducible.  An irreducible quotient of
a degenerate induced module is called a degenerate irreducible
module.

One of the most important problems of representation theory is to
determine all degenerate representations.
I will state here the result for $L=E(3|6)$
with $L_0=\Pi_{j \geq 0} \fa_j$ (see~\S\ref{sec:4}), so that the
finite-dimensional irreducible $L_0$-modules are actually
$\fa_0=s\ell_3 \oplus s\ell_2 \oplus g\ell_1$-modules (with
$\Pi_{j>0}\fa_j$ acting trivially).  We shall normalize the
generator $Y$ of $g\ell_1$ by the condition that its eigenvalue
on $\fa_{-1}$ is $-1/3$.  The finite-dimensional irreducible
$\fa_0$-modules are labeled by triples $(mn,b,Y)$, where $mn$
(resp.~$b$) are labels of the highest weight of an irreducible
representation of $s\ell_3$ (resp.~$s\ell_2$), so that $m0$ and
$0m$ label $S^m\CC^3$ and $S^m\CC^{3*}$ (resp. $b$ labels
$S^b\CC^2$), and $Y$ is the eigenvalue of the central element
$Y$. Since irreducible $E(3|6)$-modules are unique quotients of
induced modules, they can be labeled by the above triples as
well.

\begin{theorem} \cite{KR}
  \label{th:3}
  The complete list of irreducible degenerate $E(3|6)$-modules is
  as follows $(m,b \in \ZZ_+)$:
  \begin{eqnarray*}
    (0m,b,-b-\tfrac{2}{3}m-2) ,
 (0m,b,b-\tfrac{2}{3} m) ,
  (m0,b,-b+\tfrac{2}{3}m) ,
 (m0,b,b+\tfrac{2}{3}m+2)\, .
  \end{eqnarray*}
\end{theorem}

\section{Fundamental particle multiplets}
\label{sec:7}

In order to explain the connection of representation theory of
linearly compact Lie superalgebras to particle physics, let me
propose the following axiomatics of fundamental particles:

\Alphalist

\begin{enumerate}
\item %%A
  The algebra of symmetries is a linearly compact Lie
  superalgebra $L$ with an element $Y$, called the hypercharge
  operator, such that
\romanparenlistii
\begin{enumerate}
\item %%i
  $\ad Y$ is diagonalizable and normalized such that its spectrum
  is bounded below, $\subset \tfrac{1}{3} \ZZ$ and $\not\subset
  \ZZ$,

\item %%ii
  the centralizer of $Y$ in $L$ is $\fa_0=s\ell_3 + s\ell_2 +\CC
  Y$ (one may weaken this by requiring $\supset$ in place of $=$).
\end{enumerate}

\item %%B
  A particle multiplet is an irreducible subrepresentation of
  $\fa_0$ in a degenerate irreducible representation of $L$.
  Particles in a multiplet are linearly independent eigenvectors
  of the $s\ell_2$ generator $\displaystyle{I_3=\tfrac{1}{2} \left(
    \begin{array}{cccccc}
1&0\\0&-1
    \end{array}
\right)}$.  Charge $Q$ of a particle is given by the
Gell-Mann-Nishijima formula:
\begin{displaymath}
  Q=(I_3 \hbox{ eigenvalue) }+\tfrac{1}{2} \hbox{ (hypercharge).}
\end{displaymath}

\item %%C
  Fundamental particle multiplet is a particle multiplet such
  that

\romanparenlistii
\begin{enumerate}
\item %%i
 $|Q| \leq 1$ for all particles of the multiplet,

\item %%ii
only the $1$-dimensional, the two fundamental representations or the adjoint
representation of $s\ell_3$ occur.

\end{enumerate}

\end{enumerate}

Using Theorem~\ref{th:3}, it is easy to classify all fundamental
multiplets when the algebra of symmetries $L=E(3|6)$ \cite{KR}.
The answer is given in the left half of  Table~\ref{tab:1}.  The
right half contains all the fundamental particles of the
Standard model (see e.g.~\cite{O}):  the upper part is comprised
of three generations of quarks and the middle part of three
generations of leptons (these are all fundamental fermions from
which matter is built), and the lower part is comprised of
fundamental bosons (which mediate the strong and electro-weak interactions).  Except for the last line,
the match is perfect.

\begin{table}[htbp]
  \begin{center}
    \caption{}
    \label{tab:1} %%%{Table 1.}
\begin{tabular}{c c | ccc }
multiplets & charges  && particles\\
\hline \\[-1ex]
$(01,1,1/3)$ & $2/3,-1/3$ & $\binom{u_L}{d_L}$ & $\binom{c_L}{s_L}$
    & $\binom{t_L}{b_L}$\\[1ex]
$(10,1,-1/3)$ & $-2/3,1/3$ & $\binom{\tilde{u}_R}{\tilde{d}_R}$ &
     $\binom{\tilde{c}_R}{\tilde{s}_R}$ & $\binom{\tilde{t}_R}{\tilde{b}_R}$\\[1ex]
$(10,0,-4/3)$ & $-2/3$ & $\tilde{u}_L$ & $\tilde{c}_L$
    & $\tilde{t}_R$\\[1ex]
$(01,0,4/3)$ & $2/3$ & $u_R$ & $c_R$ & $ t_R$\\[1ex]
$(01,0,-2/3)$ & $-1/3$ & $d_R$ & $s_R$ &$b_R$\\[1ex]
$(10,0,2/3)$ & $1/3$ & $\tilde{d}_L$ & $\tilde{s}_L$ & $
   \tilde{b}_L$\\[-1ex]
%
%\rule{.25in}{1pt} \quad \rule{.25in}{1pt} \quad\rule{.25in}{1pt}
%\quad\rule{.25in}{1pt} \\
\setlength{\unitlength}{0.1in}
\begin{picture}(10,3)(0,-1)
  \multiput(0,0)(2,0){6}{\line(1,0){1.5}}
\end{picture}
&
\setlength{\unitlength}{0.1in}
\begin{picture}(10,3)(0,-1)
  \multiput(0,0)(2,0){6}{\line(1,0){1.5}}
\end{picture}
&
\setlength{\unitlength}{0.1in}
\begin{picture}(5,3)(1,-1)
  \multiput(0,0)(2,0){6}{\line(1,0){1.5}}
\end{picture}
&
\setlength{\unitlength}{0.1in}
\begin{picture}(5,3)(0,-1)
  \multiput(0,0)(2,0){5}{\line(1,0){1.5}}
\end{picture}
&
\setlength{\unitlength}{0.1in}
\begin{picture}(5,3)(0,-1)
  \multiput(0,0)(2,0){4}{\line(1,0){1.5}}
\end{picture}
\\
$(00,1,-1)$ & $0,-1$ & $\binom{\nu_L}{e_L}$
   & $\binom{\nu_{\mu L}}{\mu_L}$ & $\binom{\nu_{\tau
       L}}{\tau_L}$\\[1ex]
$(00,1,1)$ & $0,1$ & $\binom{\tilde{\nu}_R}{\tilde{e}_R}$
  & $\binom{\tilde{\nu}_{\mu R}}{\tilde{\mu}_R}$
  & $\binom{\tilde{\nu}_{\tau R}}{\tilde{\tau}_R}$\\[1ex]
$(00,0,2)$ & $1$ & $\tilde{e}_L$ & $\tilde{\mu}_L$
   & $\tilde{\tau}_L$\\[1ex]
$(00,0,-2)$ & $-1$ & $e_R$ & $\mu_R$ & $\tau_R$\\[1ex]
\hline \\[-1ex]
$(11,0,0)$ & $0$ & gluons\\[1ex]
$(00,2,0)$ & $1,-1,0$ & $W^+,W^-,Z$ & (gauge bosons)\\[1ex]
$(00,0,0)$ & $0$ & $\gamma$ & (photon)\\[1ex]
$(11,0, \pm 2)$ & $\pm 1$ & --

\end{tabular}
  \end{center}
\end{table}

\section{Speculations and visions}
\label{sec:8}

As the title of the conference suggests, each speaker is expected
to propose his (or her) visions in mathematics for the 21\st{st}
century.  This is an obvious invitation to be irresponsibly
speculative.  Some of the items proposed below are of this
nature, but some others are less so.

\arabiclist

\begin{enumerate}
\item %%1
  It is certainly impossible to classify all simple
  infinite-dimensional Lie algebras or superalgebras.  The most
  popular types of conditions that have emerged in the past
  30~years and that I like most are these:

\alphaparenlistii
\begin{enumerate}
\item %%a
  existence of a gradation by finite-dimensional subspaces and
  finiteness of growth \cite{K1}, \cite{M}.

\item %%b
  topological conditions \cite{G2}, \cite{K7}, \S\S~1--3,

\item %%c
  the condition of locality \cite{DK}, \cite{K4},\cite{K5}, \S~5.
\end{enumerate}

Problem (a) in the Lie algebra case has been completely solved in
\cite{M}, but an analogous conjecture in the Lie superalgebra
case \cite{KL} is apparently much harder.

Concerning (b), let me state a concrete problem.  Let $L=\CC
((x))^n$, where $\CC ((x))$ is the space of formal Laurent series
in $x$ with formal topology.  Examples of simple topological Lie
algebras with the underlying space~$L$ are the completed
(centerless) affine and Virasoro algebras.  Are there any other
examples?

Incidentally, after going to the dual, Theorem~\ref{th:1} gives a
complete classification of simple Lie co-superalgebras.

\item %%2
  In \S\ref{sec:5} I explained how to use classification of
  simple linearly compact Lie superalgebras of growth~1 in order
  to classify simple ``linear'' OPE of chiral fields in
  $2$-dimensional conformal field theory.  Will $CK_6$ play a
  role in physics or is it just an exotic animal?  Are the linearly compact
  Lie superalgebras of growth~$>1$ in any way related to OPE of
  higher dimensional quantum field theories?

\item %%3
  Each of the four types $W$, $S$, $H$, $K$ of simple primitive
  Lie algebras $(L,L_0)$ correspond to the four most important
  types of geometries of manifolds:  all manifolds, oriented
  manifolds, symplectic and contact manifolds.  Since every
  smooth supermanifold of dimension $(m|n)$ comes from a rank $n$
  vector bundle on a $m$-dimensional manifold, it is natural to
  expect that each of the simple primitive Lie superalgebras
  corresponds to one of the most important types of geometries of
  vector bundles on manifolds.  For example, the five exceptional
  superalgebras have altogether, up to conjugacy, $15$~maximal
  open subalgebras.  They correspond to irreducible
  $\ZZ$-gradations listed in \cite{CK3}:  $4$ for $E(1|6)$, $3 $
  for $E(3|6)$, $3$ for $E(3|8)$, $1$ for $E(4|4)$ and $4$ for
  $E(5|10)$ (as Shchepochkina pointed out, we missed two
  $\ZZ$-gradations of $E(1|6)$:  $(1,0,0,0,1,1,1)$ and
  $(2,2,0,1,1,1,1)$ in notation of \cite{CK3}) .  There should be
  therefore $15$~exceptional types of geometries of vector
  bundles on manifolds which are especially important.

\item %%4
  The main message of \S~\ref{sec:7} of my talk is the following
  principle:
  \begin{list}{}{}
  \item\emph{ Nature likes degenerate representations. }

  \end{list}

 There are
  several theories where this principle works very well.  First,
  it is the theory of $2$-dimensional statistical lattice models,
  especially the minimal models of \cite{BPZ}, which are (for
  $0<c<1$) nothing else but the top degenerate representations of
  the Virasoro algebra, and the WZW models, including the case of
  fractional levels, which are based on degenerate top modules
  over the affine Kac-Moody algebras (see \cite{K3} for a review
  on these modules).  Second, it is the theoretical explanation
  of the quantum Hall effect by \cite{CTZ} based on degenerate
  top modules of $W_{1+\infty}$ (see \cite{KRa}).

\item %%5
  In view of the discussion in \S\S~\ref{sec:4} and \ref{sec:7},
  it is natural to suggest that the algebra $su_3+su_2+u_1$ of
  internal symmetries of the Weinberg-Salam-Glashow Standard model extends to
  $E(3|6)$.  I am hopeful that representation theory will shed
  new light on various features of the Standard model (including
  the Kobayashi-Maskawa matrix).  It turns out \cite{KR} that all
  degenerate $E(3|6)$ Verma modules have a unique non-trivial
  singular vector.  This should lead to some canonical
  differential equations on the correlation functions (cf.~\cite{BPZ}).

I find it quite remarkable that
  the $SU_5$ Grand unified model of Georgi-Glashow combines the
  left multiplets of fundamental fermions in precisely the
  negative part of the consistent gradation of $E(5|10)$ (see
  \S~\ref{sec:4}).  This is perhaps an indication of the
  possibility that an extension from $su_5$ to $E(5|10)$ algebra
  of internal symmetries may resolve the difficulties with the
  proton decay.

One, of course, may try other finite- or infinite-dimensional Lie
superalgebras.  For example, J. van der Jeugt has tried recently
$L=s\ell(3|2)$ and it worked rather nicely, but $osp(6|2)$ has
been ruled out.

\item %%6
  Let me end with the most irresponsible suggestion.  Since $W_4$
  is, on the one hand, the algebra of symmetries of Einstein's
  gravity theory, and, on the other hand, the even part of
  $E(4|4)$, it is a natural guess that $E(4|4)$ is the algebra
  of symmetries of a nice super extension of general relativity.
  One knows that the algebra of symmetries of the minimal $N=1$
  supergravity theory is $S(4|2)$ \cite{OS}.
\end{enumerate}

\end{document}